\documentclass[11pt,a4paper]{article}

\usepackage[utf8]{inputenc}
\usepackage[T1]{fontenc}
\usepackage{lmodern}
\usepackage[english]{babel}
\usepackage{amsmath,amssymb,amsthm,mathtools}
\usepackage[a4paper,margin=27mm,top=25mm,bottom=25mm]{geometry}
\usepackage{microtype}
\usepackage{enumitem}
\usepackage[unicode,hidelinks]{hyperref}
\usepackage{authblk}
\usepackage{graphicx}

\newtheorem{theorem}{Theorem}
\newtheorem{lemma}{Lemma}
\theoremstyle{remark}
\newtheorem*{remark}{Remark}
\numberwithin{equation}{section}
\setlist[enumerate]{itemsep=0.3em,topsep=0.4em}
\setlist[itemize]{itemsep=0.25em,topsep=0.4em}
\title{Avoiding two cycles in Avoider--Enforcer games}

\author[]{Jovana Forcan}

\affil[]{Department of Computer Science and Systems, Faculty of Philosophy, \newline
University of East Sarajevo,
Bosnia and Herzegovina}

\affil[]{\texttt{jovana.forcan@ff.ues.rs.ba}}

\date{}

\begin{document}
\maketitle

\begin{abstract}
We study the strict and monotone $(1:b)$ Avoider--Enforcer games on the edges of the complete graph $K_n$. We prove that, for integers $n\ge2$ and $b\ge6n$, Avoider can keep his graph a forest after every move, except possibly the last. Thus, his final graph contains at most one cycle. This answers Question~10 of Clemens, Ehrenm\"uller, Person, and Tran, which asks whether a bias linear in $n$ is sufficient.
\end{abstract}

\section{Introduction}

Avoider--Enforcer games are positional games played on a finite
board $X$, with a family $\mathcal F\subseteq 2^X$ of \textit{losing sets}. The players, Avoider and Enforcer, alternately claim previously unclaimed elements of $X$, with Avoider moving first. Avoider loses as soon as he has claimed every element of some $F\in\mathcal F$. He wins if all elements of the board have been claimed and he has not completed any losing set.

In a Maker--Breaker game, Maker tries to claim every element of one of the sets in $\mathcal F$, while Breaker tries to prevent this. In an Avoider--Enforcer game, Avoider tries to avoid completing such a set, while Enforcer tries to force him to do so. For general background on positional games, see the books by Beck~\cite{beck2008} and by Hefetz et al.~\cite{hefetz2014}.

Early work on Avoider--Enforcer games goes back to
Lu~\cite{lu1991}. Hefetz et al.~\cite{hefetz2007} later gave a general sufficient condition for Avoider to win and applied it to several games on graphs. Further work~\cite{hefetz2008}
studied how the bias affects Avoider's ability to preserve some graph properties. 

We consider two versions of the biased $(1:b)$ game. Under the \textit{strict} rules, Avoider claims exactly one element per move, and Enforcer claims exactly $b$. Under the \textit{monotone} rules, they claim at least one and at least $b$ elements, respectively. If fewer than $b$ elements remain on Enforcer's turn, he claims all of them. Under the strict rules, increasing the bias does not necessarily
favor Avoider. To ensure monotonicity in the bias, Hefetz et al.~\cite{hefetz2010} introduced the monotone rules.

In this paper, the board is $E(K_n)$, and Avoider's graph consists of his claimed edges together with all $n$ vertices. We aim to keep this graph as close to a forest as possible. Clemens et al.~\cite{clemens2015} proved that, for sufficiently large $n$ and $b\ge 200n\ln n$,
Avoider can keep his graph a forest after every move except possibly the last, under both sets of rules. His final graph then contains at most one cycle. In Question~10 of the same paper, they asked whether a bias linear in $n$ is sufficient
to obtain at most one cycle.

We answer this question affirmatively for both versions of the game.

\begin{theorem}\label{thm:main}
Let $n\ge2$ and $b\ge6n$ be integers. In both the strict and the monotone $(1:b)$ Avoider--Enforcer game on $E(K_n)$, Avoider has a strategy such that his graph is a forest after
every one of his moves except possibly the last. His final graph is a forest together with at most one additional edge and therefore contains at most one cycle in total.
\end{theorem}

\subsection{Notation}
\label{subsec:notation}

Write $E_A$ for the set of edges already claimed by Avoider. Throughout the game, his graph is the \textit{spanning} graph
\begin{equation}
  A=(V,E_A).
  \label{eq:spanning}
\end{equation}
The players claim edges, while the vertex set $V$ stays fixed.
Thus $A$ contains all $n$ vertices throughout the game. Vertices not incident with any of Avoider's edges are isolated, and each forms a component of size one.

Let $F\subseteq E(K_n)$ be the current set of free edges. The pair $(A, F)$ describes the board at that point in the game.
If $E_E$ is the set of edges claimed by Enforcer, then
$$
  E(K_n)=E_A\,\dot\cup\,E_E\,\dot\cup\,F,
  \qquad
  E_E=E(K_n)\setminus(E_A\cup F).
$$
Thus $A$ and $F$ also determine Enforcer's edges. 
At the start of the game, $E_A=\varnothing$ and $F=E(K_n)$.

For the quantities $B, \Phi, \Psi, I$ that will be defined later, a symbol without a subscript denotes its value immediately before Avoider's move.
The subscript $A$ denotes its value immediately after that move, before Enforcer's reply. The subscript $+$ denotes its value immediately after Enforcer's reply. If the game continues, this is also the value immediately before Avoider's next move.

The subscript $0$ denotes the value at the start of the game, before either player has claimed an edge.

\section{Component structure}
\label{sec:positions}

Let $T_1,\ldots,T_m$ be the connected components of the
current spanning graph $A=(V,E_A)$, and put
$$
  s_i=|V(T_i)|.
$$
Every vertex belongs to exactly one component, so
\begin{equation}
  V=V(T_1)\,\dot\cup\,\cdots\,\dot\cup\,V(T_m),
  \qquad
  \sum_{i=1}^{m}s_i=n.
  \label{eq:partition}
\end{equation}
This equality holds throughout the game, since isolated vertices are also counted as components of size one.

A free edge is called \textit{external} if its endpoints lie in different components of $A$, and \textit{internal} if its endpoints lie in the same component of $A$.

For distinct $i,j$, let
$$
  f_{ij}
  =\bigl|\{uv\in F:\ u\in V(T_i),\ v\in V(T_j)\}\bigr|.
$$
In particular, $f_{ji}=f_{ij}$. There are exactly $s_i s_j$
edges of the original complete graph between vertex sets $T_i$ and $T_j$, of which only $f_{ij}$ are free. Hence
\begin{equation}
  0\le f_{ij}\le s_i s_j
  \qquad(i\ne j).
  \label{eq:capacity}
\end{equation}
This inequality holds throughout the game for the components at that point. When components merge, their sizes and the numbers $f_{ij}$ change, but the inequality still holds.

\subsection{Avoider's strategy}
\label{subsec:strategy}

Avoider always claims exactly one edge. His strategy is:
\begin{enumerate}[label=\textnormal{(\roman*)}]
  \item If a free external edge exists, choose distinct components $T_i,T_j$ with $f_{ij}>0$ for which $s_i s_j$ is minimum, and claim any free edge between them.
  \item If no free external edge exists but some edge is still free, claim any free internal edge.
\end{enumerate}
If several pairs in \textnormal{(i)} have the same minimum product, Avoider may choose any of them. If no free edge remains, the game has already ended.

If $A$ is a forest, a move of type \textnormal{(i)} joins two trees and keeps the graph a forest. We will prove that a move of type \textnormal{(ii)}, if it occurs, is Avoider's last move.

\begin{figure}[!h]
\centering
\includegraphics[width=0.95\textwidth]{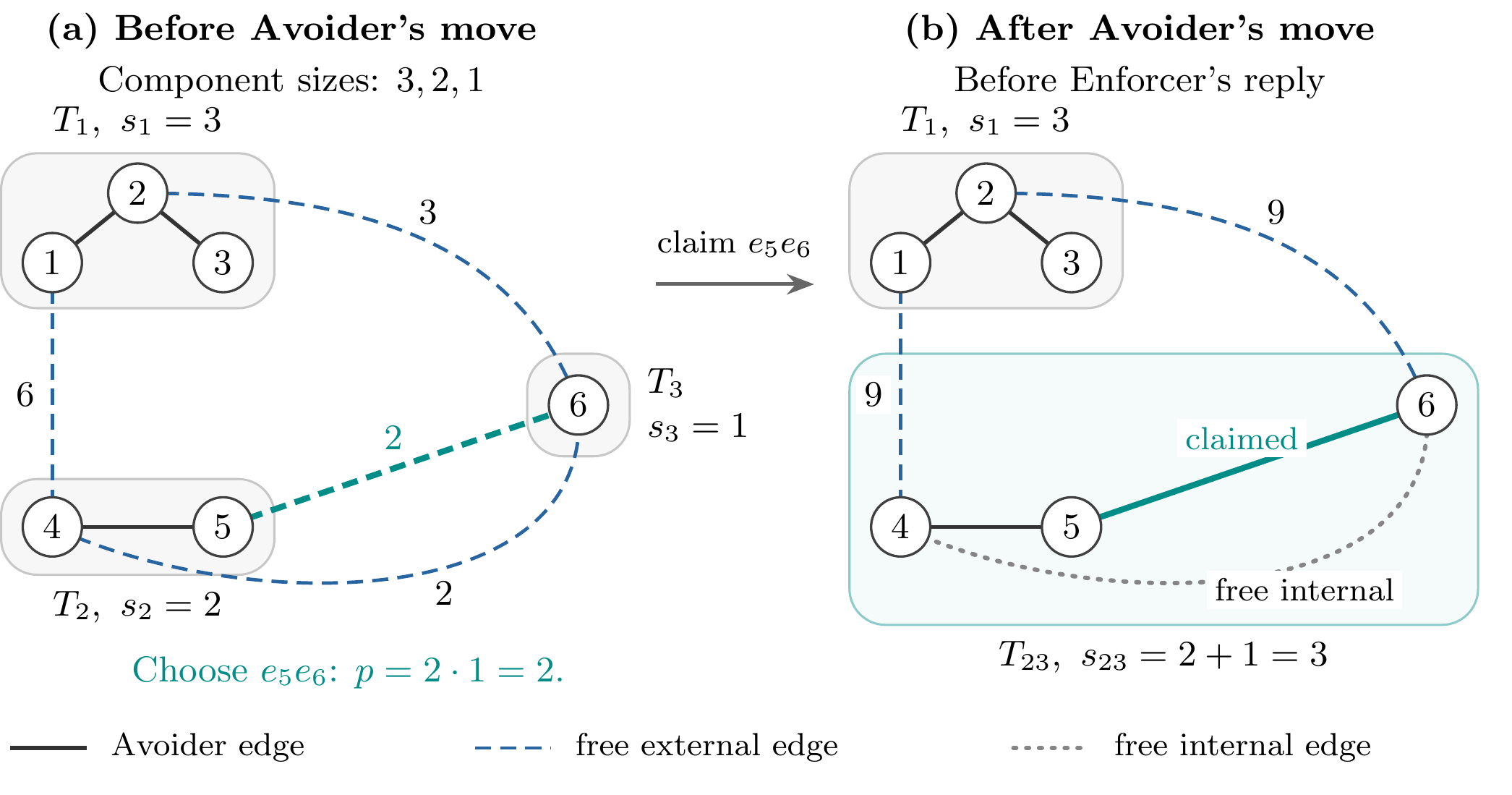}
\caption{An example of Avoider's strategy.}
\label{fig:strategy}
\end{figure}

Figure~\ref{fig:strategy} illustrates a position with component sizes $3,2,1$. All three pairs of components have a free edge, and their products are $6,3,2$.
Both $e_4e_6$ and $e_5e_6$ have the minimum product $p=2$. Avoider may choose either. The figure~\ref{fig:strategy} shows the choice of $e_5e_6$. The two remaining components have sizes $3$ and $3$,
so any free edge between them now has product $9$. The edge $e_4e_6$ remains free but becomes internal.

\begin{remark}[Contraction]
We can view each component $T_i$ as a single vertex labelled by its size $s_i$. We keep all $f_{ij}$ free original edges between each pair of components. After contraction, these may be parallel edges between the two new vertices.

When $T_i$ and $T_j$ merge, the new component has size
$s_i+s_j$. The number of free edges between it and any other component $T_k$ is $f_{ik}+f_{jk}$. We count free internal edges separately. All edge counts refer to original edges of $K_n$. Both $n$ and the bias $b$ stay fixed throughout the game.
\end{remark}

\section{The potential function}
\label{sec:potential}

\subsection{Definition}
To analyze Avoider's strategy, we use a potential function based on the component sizes and the free edges between them.

For a free external edge $uv$ between $T_i$ and $T_j$, define its weight by
\begin{equation}
w(uv)=s_i s_j.
\label{eq:edgeweight}
\end{equation}
It depends on the current component sizes and can change when components merge.
Every free edge between the same pair of components has the same weight.

For example, between components of sizes $3$ and $2$, each free edge has weight $6$. If Enforcer has claimed two of the six original edges between these components, then $f_{ij}=4$.
The product $s_i s_j$ is still $6$, and the four free edges contribute $4\cdot6=24$ to the sum of weights.

Define
\begin{equation}
B=\sum_{i=1}^{m}s_i^2,\qquad
\Phi=\sum_{1\le i<j\le m}f_{ij}s_i s_j,\qquad
\Psi=\Phi+nB.
\label{eq:potentials}
\end{equation}
Thus, $\Phi$ is the sum of the weights of all free external edges. Free internal edges do not contribute to $\Phi$.

Although $\sum_i s_i=n$ is fixed, $B=\sum_i s_i^2$ is not.
Initially $B=n$. After Avoider claims his first edge, there is one component of size two and $n-2$ isolated vertices, giving
$$
B=2^2+(n-2)\cdot1^2=n+2.
$$
In general, merging components of sizes $x,y$ increases $B$
by $(x+y)^2-x^2-y^2=2xy$. Enforcer's moves do not change Avoider's components and therefore do not change $B$.

We use $\Psi=\Phi+nB$ to bound $B$ and, through $B$, the number of free internal edges (see Lemma~\ref{lem:internal}).
When Avoider joins two components, $B$ increases.
Enforcer does not change $B$, but he decreases $\Phi$ by claiming free external edges.
In Section~\ref{sec:round}, we show that, under the stated assumptions, the decrease in $\Phi$ over a full round is at least the increase in $nB$. Therefore, $\Psi$ does not increase over that round.

Here $\Psi$ denotes the value immediately before Avoider's move, and $\Psi_+$ denotes the value immediately after Enforcer's reply. The inequality,
$$
\Psi_+\le\Psi,
$$
compares the beginning and the end of a full round. The value of $\Psi$ may increase immediately after Avoider's move.

\subsection{The invariant}

Initially every component has size one and all edges are free, so
\begin{equation}
B_0=n,\qquad
\Phi_0=\binom n2,\qquad
\Psi_0=\binom n2+n^2=\frac{3n^2-n}{2}.
\label{eq:initial}
\end{equation}

We will prove by induction that
\begin{equation}
\Psi\le\Psi_0
\label{eq:invariant}
\end{equation}
immediately before each Avoider move, while his graph is a forest.

At the start of the game, \eqref{eq:invariant} holds with equality. In Section~\ref{sec:round}, we prove that a full round preserves this inequality: if it holds before Avoider's move, then it also holds after Enforcer's reply.
The induction is completed in the proof of
Theorem~\ref{thm:main}, where we also consider the cases in which the game ends.

\section{Merging two components}
\label{sec:merge}

\begin{lemma}
\label{lem:internal}
At any position $(A,F)$, let $I$ be the number of free
internal edges. Then
\begin{equation}
I\le\sum_{i=1}^{m}\binom{s_i}{2}
=\frac{B-n}{2}.
\label{eq:internal}
\end{equation}
\end{lemma}

\begin{proof}
There are at most $\binom{s_i}{2}$ edges within the vertex set of $T_i$, and its free internal edges form a subset of these. Summing over all components and using \eqref{eq:partition},
$$
  I\le\frac12\sum_i(s_i^2-s_i)=\frac{B-n}{2}.
$$
This upper bound may also count edges that have already been claimed, which only makes the upper bound larger.
\end{proof}

\begin{lemma}
\label{lem:merge}
Suppose that Avoider claims a free edge between $T_i,T_j$ according to the minimum product rule. Write
$$
x=s_i,\qquad y=s_j,\qquad p=xy,\qquad f=f_{ij}\ge1.
$$
Then
\begin{align}
B_A&=B+2p,
\label{eq:Bmerge}\\
\Phi_A-\Phi
 &=y\sum_{k\ne i,j}f_{ik}s_k
+x\sum_{k\ne i,j}f_{jk}s_k-fp,
 \label{eq:exact}\\
\Phi_A-\Phi
&\le2p(B-x^2-y^2)-fp
\le2pB.
\label{eq:phibound}
\end{align}
Moreover, every free external edge after Avoider's move has weight at least $p$.
\end{lemma}

\begin{proof}
The two components are replaced by a component of size $x+y$, while all other component sizes remain unchanged. Therefore
$$
B_A-B=(x+y)^2-x^2-y^2=2xy=2p,
$$
proving \eqref{eq:Bmerge}.

To derive \eqref{eq:exact}, first consider the $f$ free edges between $T_i$ and $T_j$ before Avoider's move. Their total contribution to $\Phi$ is $fxy=fp$. Avoider claims one of
them, and the other $f-1$ become internal edges of the new component. None of these edges contributes to $\Phi_A$.
Their total contribution to the change is $-fp$. The remaining $f-1$ edges are still free, but they no longer contribute to $\Phi_A$ because they are now internal.

For each $k\ne i,j$, there are now $f_{ik}+f_{jk}$ free edges between the merged component and $T_k$. These edges change the potential by
\begin{align*}
 &(f_{ik}+f_{jk})(x+y)s_k
       -f_{ik}xs_k-f_{jk}ys_k\\
 &\hspace{2em}=y f_{ik}s_k+x f_{jk}s_k.
\end{align*}
Edges between components that do not merge contribute the same amount as before. Previously internal edges remain internal.
Adding these contributions proves \eqref{eq:exact}.

For the first inequality in \eqref{eq:phibound}, use
$$
f_{ik}\le xs_k,\qquad f_{jk}\le ys_k
\qquad(k\ne i,j).
$$
It follows that
\begin{align*}
 y\sum_{k\ne i,j}f_{ik}s_k
 +x\sum_{k\ne i,j}f_{jk}s_k
 &\le y\sum_{k\ne i,j}xs_k^2
       +x\sum_{k\ne i,j}ys_k^2\\
 &=2xy\sum_{k\ne i,j}s_k^2\\
 &=2p(B-x^2-y^2).
\end{align*}
The last equality uses $\sum_{k\ne i,j}s_k^2=B-x^2-y^2$, because we leave out the two components being joined.
Keeping the term $-fp$ from \eqref{eq:exact} gives the first
inequality. For the second inequality, note that
$$
2p(B-x^2-y^2)-fp
=2pB-\bigl(2p(x^2+y^2)+fp\bigr)\le2pB,
$$
since the subtracted quantity is nonnegative.

Finally, consider any edge that is still free and external after Avoider's move. Before the move, it was also free and external. Its weight was therefore at least $p$, since Avoider chose
an edge of minimum weight.

When two components merge, the sizes of the components containing the endpoints of this edge either stay the same or increase. Since the weight is the product of these sizes, it cannot decrease. Thus the edge still has weight at least $p$.

The remaining free edges between the two merged components are now internal, so the statement about external edges does not apply to them. This proves the lemma.
\end{proof}

\section{Preserving the invariant over a full round}
\label{sec:round}

A \textit{full round} in the following argument consists of an Avoider move joining two components, followed by an Enforcer reply in which at least $b$ edges are claimed.

\begin{lemma}
\label{lem:round}
Let $b\ge6n$. Suppose that immediately before Avoider's move his graph is a forest and $\Psi\le\Psi_0$. If he joins two components according to the strategy and Enforcer then claims at least $b$ edges, then
\begin{equation}
  \Psi_+\le\Psi\le\Psi_0.
  \label{eq:round}
\end{equation}
\end{lemma}

\begin{proof}
Use the notation $x,y,p$ of Lemma~\ref{lem:merge}.
By the assumption $\Psi\le\Psi_0$ and since $\Phi\ge0$,
we have
$$
  nB\le\Phi+nB=\Psi\le\Psi_0.
$$
Dividing by $n>0$, we obtain
\begin{equation}
B\le\frac{\Psi_0}{n}
=\frac{3n-1}{2}<\frac{3n}{2}.
\label{eq:Bbound}
\end{equation}
Also,
$$
  2p=2xy\le x^2+y^2\le B.
$$

Together with \eqref{eq:Bmerge}, this gives the  bound
\begin{equation}
B_A=B+2p\le2B<3n.
\label{eq:temporary}
\end{equation}

Apply Lemma~\ref{lem:internal} after Avoider's move, using $I_A$ and $B_A$. We get
\begin{equation}
I_A\le\frac{B_A-n}{2}<n.
\label{eq:IA}
\end{equation}
This includes all remaining edges that have just become internal as a result of the merging.

Suppose Enforcer claims $d\ge b\ge6n$ edges, of which $q$ are external. At most $I_A$ can be internal. So,
$$
  q\ge d-I_A\ge b-I_A>5n,
$$
and in particular $q\ge5n$.

Enforcer's move does not change Avoider's components. Therefore, the weights of free external edges are fixed during that reply.
Every external edge Enforcer claims has weight at least $p$, by Lemma~\ref{lem:merge}, and its removal subtracts that weight from $\Phi_A$. Claiming an internal edge does not change $\Phi$. Consequently,
\begin{equation}
\Phi_+\le\Phi_A-qp\le\Phi_A-5np.
\label{eq:Edecrease}
\end{equation}

On the other hand, \eqref{eq:phibound} and \eqref{eq:Bbound} give
$$
  \Phi_A-\Phi\le2pB\le3np.
$$
Substituting in \eqref{eq:Edecrease}, we obtain
\begin{equation}
\Phi_+\le\Phi-2np.
  \label{eq:phiround}
\end{equation}
Since $B_+=B_A=B+2p$, it follows that
\begin{align*}
\Psi_+
&=\Phi_++nB_+\\
&\le(\Phi-2np)+n(B+2p)\\
&=\Phi+nB=\Psi\le\Psi_0.
\end{align*}
This proves the lemma.
\end{proof}

\begin{remark}
The proof uses $B<3n/2$ before Avoider's move.
After that move, the weaker bound $B_A<3n$ is enough to obtain $I_A<n$.
\end{remark}

\section{Proof of the theorem}
\label{sec:finish}

\begin{proof}[Proof of Theorem~\ref{thm:main}]
Avoider uses the strategy in Section~\ref{subsec:strategy}.
Initially, his graph is edgeless, hence a forest, and
$\Psi=\Psi_0$ by \eqref{eq:initial}.

We proceed inductively through the positions immediately before Avoider's moves. Assume that at such a position $A$
is a forest and $\Psi\le\Psi_0$. If no edges remain, the game
has already ended. Otherwise, there are two cases.

\medskip
\noindent\textbf{Case 1: a free external edge exists.}
Avoider claims a free edge between two components with the minimum product of their sizes. Joining two trees by one edge produces a tree, so his graph remains a forest.

If no free edge remains after this move, the game ends
immediately with a forest. If between one and $b-1$ free
edges remain, Enforcer must claim all of them, again ending
the game with a forest.

If at least $b$ free edges remain, Enforcer claims at least
$b$ edges. Lemma~\ref{lem:round} gives
$\Psi_+\le\Psi_0$. If Enforcer has claimed all remaining edges, the game ends with a forest. Otherwise, Avoider is to move again. His graph is still a forest and the invariant still holds, so the induction can continue.

\medskip
\noindent\textbf{Case 2: no free external edge exists.}
All remaining free edges are internal. Let their number be
$I\ge1$. Lemma~\ref{lem:internal} and \eqref{eq:Bbound} give
\begin{equation}
I\le\frac{B-n}{2}
\le\frac{n-1}{4}<b.
\label{eq:terminal}
\end{equation}
Avoider claims one of these edges, say $uv$.

The endpoints $u,v$ belong to the same tree of the current forest. There is a unique path between them in that tree. Adding the free edge $uv$ creates exactly one cycle, consisting of this path and the new edge. All other components
remain trees.

If $I=1$, no free edges remain after Avoider's
move. If $I>1$, there are $I-1<b$ free edges left, and Enforcer must claim all of them. In either case, Avoider never moves again and cannot create a second cycle.

Thus Avoider's graph remains a forest after every move except possibly his last, and that last move can add only one edge to a forest.

The same proof works for the strict and monotone games. Choosing exactly one edge is legal for Avoider in both versions. For every full Enforcer reply, we used only $d\ge b$. In both versions, Enforcer must claim all remaining edges if fewer than $b$ are left. Thus the result holds in both games.
\end{proof}

\appendix
\section{A sharper potential estimate}
\label{sec:sharper}

The estimate in this appendix is not needed for the proof of Theorem~\ref{thm:main}. A more detailed calculation gives a better bound on the decrease of the potential over a full round.

Assume the hypotheses of Lemma~\ref{lem:round}, and use
$$
p=xy>0,\qquad f=f_{ij}\ge1.
$$
As before, quantities without a subscript refer to the position before Avoider's move.

Lemma~\ref{lem:internal}, applied after Avoider's move, gives
$$
I_A\le\frac{B_A-n}{2}.
$$
Substituting $B_A=B+2p$ from \eqref{eq:Bmerge}, we obtain
\begin{equation}
I_A\le\frac{B+2p-n}{2}.
\label{eq:IAexact}
\end{equation}
The total number of original vertices is still $n$ after the merging. Only the sum of squared component sizes has
changed from $B$ to $B_A$.

We calculate the change in $\Psi$ by considering each player's effect on $\Phi$ and the change in $nB$:
\begin{equation}
\begin{aligned}
\Psi_+-\Psi
&=(\Phi_++nB_+)-(\Phi+nB)\\
&=(\Phi_+-\Phi)+n(B_+-B) \\
& =(\Phi_A-\Phi)+(\Phi_+-\Phi_A)+n(B_+-B).
\end{aligned}
\label{eq:split}
\end{equation}
For the first contribution, \eqref{eq:phibound} gives
$$
\Phi_A-\Phi\le2p(B-x^2-y^2)-fp.
$$
For the second, at least $b-I_A$ of Enforcer's claimed edges are external, and each has weight at least $p$.
So,
$$
\Phi_+-\Phi_A\le-p(b-I_A).
$$
For the third, $B_+=B_A=B+2p$, so
$$
n(B_+-B)=2np.
$$
Putting these estimates into \eqref{eq:split} gives
\begin{equation}
\Psi_+-\Psi
\le2p(B-x^2-y^2)-fp-p(b-I_A)+2np.
\label{eq:sharperraw}
\end{equation}

Since $-p(b-I_A)=-pb+pI_A$, the coefficient of $I_A$
is positive. We can use the bound on $I_A$ from \eqref{eq:IAexact}. First, we take $p$ outside the brackets:
\begin{align}
 \Psi_+-\Psi
&\le p\bigl[2B-2(x^2+y^2)-f-b+I_A+2n\bigr]\notag\\
&\le p\left[
2B-2(x^2+y^2)-f-b
+\frac{B+2p-n}{2}+2n
\right]\notag\\
&=p\left[
\frac52B+\frac32n-b+p-2(x^2+y^2)-f
\right].
\label{eq:sharper}
\end{align}

Using $p=xy$, write
$$
x^2+y^2=(x-y)^2+2p.
$$
It follows that
\begin{align*}
p-2(x^2+y^2)-f
&=p-2\bigl((x-y)^2+2p\bigr)-f\\
&=-2(x-y)^2-3p-f\\
&\le0.
\end{align*}
Removing this nonpositive term makes the upper bound larger,
so \eqref{eq:sharper} implies
\begin{equation}
  \Psi_+-\Psi
  \le p\left(\frac52B+\frac32n-b\right).
  \label{eq:sharperclean}
\end{equation}

Finally, use $B\le(3n-1)/2$ and $b\ge6n$, hence
$-b\le-6n$, to obtain
\begin{align*}
\Psi_+-\Psi
&\le p\left[
    \frac52\cdot\frac{3n-1}{2}
    +\frac32n-6n
    \right]\\
&=p\left(\frac{15n-5+6n-24n}{4}\right)\\
&=-\frac{3n+5}{4}\,p<0.
\end{align*}
The last inequality is strict because $p>0$. Thus every full round under the hypotheses of Lemma~\ref{lem:round} decreases $\Psi$ by at least $(3n+5)p/4$. 

\section*{Acknowledgments}

The author thanks Mirjana Mikalački for helpful discussions about the open question considered in this paper. \\
ChatGPT (OpenAI) was used to assist with developing parts of the proof and improving the presentation of the paper.
The author takes full responsibility for all proofs, calculations, and conclusions.

\end{document}